\documentclass{iciam}

\usepackage{graphicx}    
\usepackage{color}

\usepackage{tikz}

\contact[jmsanzserna@gmail.com\ ]{Departamento de Matem\'aticas, Universidad Carlos III de Madrid, E--28911 Legan\'es (Madrid),  Spain.}

\contact[Ander.Murua@ehu.es]{Konputazio Zientziak eta A.\ A.\  Saila, Informatika
 Fakultatea, UPV/EHU, E--20018 Donostia--San Sebasti\'{a}n,  Spain.}

\newcommand{\C}{\mathbb{C}}
\newcommand{\calC}{\mathcal{C}}
\newcommand{\bk}{{\bf k}}
\newcommand{\trees}{\mathcal{T}}
\newcommand{\eldiff}{\mathcal{F}}
\newcommand{\g}{\mathfrak{g}}
\newcommand{\G}{\mathcal{G}}
\newcommand{\Ham}{\mathcal{H}}
\newcommand{\R}{\mathbb{R}}
\newcommand{\W}{\mathcal{W}}
\newcommand{\Z}{\mathbb{Z}}
\newcommand{\uno}{{\, 1\!\! 1\,}}
\newcommand{\sh}{{\,\scriptstyle \sqcup\!\sqcup\,}}

\tikzstyle dtree=[grow'=up,sibling distance=4mm,level distance=4mm,thick]
\tikzstyle dtree node=[scale=0.6,shape=circle,very thin,draw]
\tikzstyle dtree black node=[style=dtree node,fill=black]
\tikzstyle dtree white node=[style=dtree node,fill=white]
\newcommand{\rta}{
\begin{tikzpicture}[dtree]
    \node[dtree black node] {}
    ;
  \end{tikzpicture}
}
\newcommand{\rtb}{
 \begin{tikzpicture}[dtree]
    \node[dtree black node] {}
   child { node[dtree black node]{} }
    ;
  \end{tikzpicture}
}
\newcommand{\rtc}{
\begin{tikzpicture}[dtree]
    \node[dtree black node] {}
    child { node[dtree black node]{}
            child { node[dtree black node]{} }
          }
    ;
  \end{tikzpicture}
  }
\newcommand{\rtd}{
 \begin{tikzpicture}[dtree]
    \node[dtree black node] {}
    child { node[dtree black node]{} }
    child { node[dtree black node]{} }
    ;
  \end{tikzpicture}
 }
\newcommand{\rte}{
    \begin{tikzpicture}[dtree]
    \node[dtree black node] {}
    child { node[dtree black node]{}
            child { node[dtree black node]{}
                   child { node[dtree black node]{} }
                  }
          }
    ;
  \end{tikzpicture}
 }
\newcommand{\rtf}{
 \begin{tikzpicture}[dtree]
    \node[dtree black node] {}
    child { node[dtree black node]{}
    child { node[dtree black node]{} }
    child { node[dtree black node]{} }
    }
    ;
  \end{tikzpicture}
  }
\newcommand{\rtg}{
 \begin{tikzpicture}[dtree]
    \node[dtree black node] {}
    child { node[dtree black node]{} }
    child { node[dtree black node]{}
            child { node[dtree black node]{} }
    }
    ;
  \end{tikzpicture}
                 }
 \newcommand{\rth}{
 \begin{tikzpicture}[dtree]
    \node[dtree black node] {}
    child { node[dtree black node]{} }
    child { node[dtree black node]{} }
    child { node[dtree black node]{} }
    ;
  \end{tikzpicture}
 }
\newlength{\treewidth}
\newlength{\treeheight}

\newcommand{\tb}[2][0.01]{
       \settowidth{\treewidth}{#2}%
        \settoheight{\treeheight}{#2}%
        \raisebox{-#1%
        \treeheight}{\mbox{\hspace{-0.1\treewidth}\rule{0. \treewidth}{1.1 \treeheight} #2 \hspace{-0.1\treewidth}}}}

\newcommand{\ts}[2][0.1]{
       \settowidth{\treewidth}{#2}%
        \settoheight{\treeheight}{#2}%
        \raisebox{-#1%
        \treeheight}{\resizebox{0.3 \treewidth}{!}{\mbox{\hspace{-0.15\treewidth} #2 \hspace{-0.15\treewidth}}}}}

\title[Formal series and numerical integrators]{Formal series and numerical integrators:  some history and  some new techniques}

\author[Sanz-Serna and Murua]{Jes\'us Mar\'{\i}a Sanz-Serna and Ander Murua}

\begin{document}

\begin{abstract} This paper provides a brief history of B-series and the associated Butcher group and presents the new theory of word series and extended word series.
B-series (Hairer and Wanner 1976) are formal series of functions parameterized by rooted trees. They  greatly simplify the study of Runge-Kutta schemes and other numerical integrators. We examine the problems that led to the introduction of B-series and survey a number of more recent developments, including applications outside  numerical mathematics. Word series (series of functions parameterized by words from an alphabet) provide in some cases a very convenient alternative to B-series. Associated with word series is a group $\G$ of coefficients with a composition rule simpler than the corresponding rule in the Butcher group.
From a more mathematical point of view, integrators, like Runge-Kutta schemes, that are affine equivariant are represented by elements of the Butcher group,  integrators that are equivariant with respect to arbitrary changes of variables are represented by elements of the word group $\G$..
\end{abstract}

\begin{classification}
Primary 65L05; Secondary 34C29, 70H05, 16T05.
\end{classification}

\begin{keywords}
B-series, word series, extended word series,  Lie groups, Lie algebras, Hamiltonian problems,  averaging, normal forms, splitting algorithms, processing numerical methods, modified systems, oscillatory problems, resonances.
\end{keywords}

\maketitle

\section{Introduction}

This contribution presents a survey of the use of formal series in the analysis of numerical integrators. The first part of the paper (Section 2) reviews a number of developments (the study of the order of consistency of Runge-Kutta  methods, Butcher's notions of composition of Runge-Kutta schemes and effective order, etc.) that underlie the introduction by Hairer and Wanner \cite{HW} of the concept of {\em B-series} in 1976. B-series are formal series parameterized by rooted trees. We also show how B-series and the associated Butcher group have gained prominence in view not only of their relevance to the analysis of geometric integrators but also of their applications to several mathematical theories not directly related to numerical analysis. In the second part of the paper (Sections 3--5) we summarize the recent theory of word series and extended word series developed in \cite{words}. These series are parameterized by words from an alphabet, rather than by rooted trees, and may be used advantageously to analyze some numerical integrators, including splitting algorithms. {\em Word series} appeared implicitly \cite{part2} and explicitly in \cite{orlando}, \cite{part3}.  {\em Extended word series} were  introduced in \cite{words}.  While series of {\em differential operators}  parameterized by words (Chen-Fliess series) are very
common in control theory  and dynamical systems (see e.g.\ \cite{nuevo}, \cite{fm}) and have also been used in numerical analysis (see, among others, \cite{tesis}, \cite{anderfocm}, \cite{k}), word series, as defined here, are, like B-series, series of {\em functions}.

Word series and extended word series  are, when applicable, more convenient than B-series. Associated with word series is a group $\G$ of coefficients with a composition rule simpler than the corresponding rule in the Butcher group.
  From a more mathematical point of view, integrators, like Runge-Kutta schemes, that are affine equivariant are represented by elements of the Butcher group,  integrators that are equivariant with respect to arbitrary changes of variables are represented by elements of the word group $\G$.

It should also be mentioned that  B-series and word series have been recently applied to reduce efficiently dynamical systems to normal forms, to perform high-order averaging and to find formal conserved quantities
\cite{part1}, \cite{part2}, \cite{orlando}, \cite{part3}, \cite{words}. The  approach via B- and word series also makes it possible to simplify the derivation of the corresponding error estimates.
These developments provide  applications outside numerical mathematics of tools originally devised to analyze numerical methods.

\section{Series based on rooted trees}
\label{s:trees}

This section contains an overview of the use B-series in the analysis of numerical integrators.

\subsection{Finding the order of a Runge-Kutta method}

A Runge-Kutta (RK) method with $s$ stages is specified by $s+s^2$ real coefficients $b_i$,  $a_{ij}$, $i,j = 1,\dots,s$. When the method is applied with stepsize $h$ to  the initial value problem\footnote{Our attention is restricted to deterministic differential equations in Euclidean spaces. We do not consider the extension to  stochastic differential equations, differential equations on Lie groups, differential-algebraic equations, etc.}
\begin{equation}\label{eq:ode}
\frac{d}{dt} x = f(x),\qquad x(0) = x_0\in\R^D,
\end{equation}
and the approximation $x_n$ to the exact solution value $x(nh)$, $n = 0, 1,\dots$, has been computed,
the formulas to find the next approximation $x_{n+1}$ are
\begin{equation}\begin{split}
X_{n,i} &= x_n + h \sum_{j=1}^s a_{ij} f(X_{n,j}), \qquad 1\leq i\leq s,\\
x_{n+1} & = x_n + h \sum_{i=1}^s b_i f(X_{n,i})
\end{split}
\label{eq:rk}
\end{equation}
(the $X_{n,i}$ are the auxiliary internal stages). If $f$ is well behaved and $|h|$ is small, the relations (\ref{eq:rk}) define $x_{n+1}$ as a function of $x_n$, i.e.\ $x_{n+1} =\psi_h(x_n)$. The $f$-dependent mapping $\psi_h$  is an approximation to the solution flow $\phi_h$ of (\ref{eq:ode}) for which
$x((n+1)h) = \phi_h(x(nh))$;  for smooth $f$ and an RK method of order $\nu$, $\psi_h(x)-\phi_h(x) = \mathcal{O}(h^{\nu+1})$ as $h\rightarrow 0$. To determine the value of $\nu$ for a given method it is therefore necessary to write down the expansion of $\psi_h(x)$ in powers of $h$, a job  trickier than one may think: it took mathematicians more than fifty years \cite{buthis} to come up with an easy, systematic way of carrying it out.
After the  work of  Butcher \cite{butcher63} ---following on earlier developments by Gill and Merson--- the Taylor series for $\psi_h$ is written, with the help of rooted trees, in the form
\begin{equation}
\psi_h(x) = x + \sum_{n=1}^\infty h^n\sum_{u\in\trees_n} \frac{1}{\sigma(u)}\, c_u\, \eldiff_u(x).
\label{eq:psi}
\end{equation}

\begin{table}[t!]
\centering
\caption{Rooted trees with $\leq 4$ vertices. For each tree the root is the bottom vertex}
\label{tab:1}       
\begin{center}
$\renewcommand\arraystretch{1.5}
   \begin{array}{c|cccccccccccccccccc|}
    |u| & 1 & 2 & 3 & 3 & 4 & 4 & 4 & 4 \\ \hline
    \sigma(u) & 1 & 1& 1& 2 & 1&  2 & 1 & 6\\ \hline
    u! & 1 & 2 & 6 & 3 & 24 & 12 & 8 & 4
    \\
    \hline
     & \rta & \rtb &  \rtc & \tb{\rtd}  & \tb{\rte} & \tb{\rtf} & \tb{\rtg} & \tb{\rth} \\ \hline
 \end{array}
$
\end{center}
\end{table}

Here:
\begin{itemize}
\item $\trees_n$ denotes the set of all rooted trees with $n$ vertices (see Table~\ref{tab:1}) and, for each rooted tree $u$,
$\sigma(u)$ is the cardinal of the group of symmetries of $u$.

\item $\eldiff_u$ (the elementary differential associated with $u$) is an $\R^D$-valued function of $x$ that depends on (\ref{eq:ode})
but does not change with the RK coefficients $b_i$, $a_{ij}$. Using the rooted tree $\ts[0.0]{\rtg}$ as an example,
\begin{equation*}
\eldiff_{\ts[0.2]{\rtg}}(x) = \partial_{xx} f(x) \Big[f(x),\partial_xf(x) \big[f(x)\big]\Big],
\end{equation*}
where $\partial_{xx}f(x)[\cdot,\cdot]$ and $\partial_xf(x) [\cdot]$ respectively denote the second and first Frechet derivatives of $f$ evaluated at $x$. The key point is to observe
how the structure of the elementary differential mimics that of the corresponding rooted tree.

\item $c_u$ (the elementary weight associated with $u$) is a a real number that changes with  the coefficients $b_i$, $a_{ij}$ and is independent of the system
(\ref{eq:ode}) being integrated. Taking again $\ts[0]{\rtg}$ as an example,
\begin{equation*}
c_{\ts[0.2]{\rtg}} = \sum_{i=1}^s b_i \Big(\sum_{j=1}^s a_{ij}\:\sum_{k=1}^s a_{ik}\big(\sum_{\ell=1}^s a_{k\ell}\big)\Big);
\end{equation*}
 the structure of the rooted tree is reflected in the summations.
\end{itemize}

For the solution flow $\phi_h$ of (\ref{eq:ode}), the expansion is
\begin{equation}\label{eq:phi}
\phi_h(x) = x + \sum_{n=1}^\infty h^n\sum_{u\in\trees_n} \frac{1}{\sigma(u)}\, \frac{1}{u!}\, \eldiff_u(x).
\end{equation}
where $u!$ (the density of $u$) is an integer that is easily computed recursively.
By comparing (\ref{eq:psi}) and (\ref{eq:phi}) we conclude that an RK method has order $\geq \nu$ if and only if $c_u =1/u! $ for each $u\in\trees_n$, $n\leq \nu$. As an illustration, we note that for order $\geq 3$ the coefficients have to satisfy the following set of {\em order conditions}
\begin{equation*}
\sum_i b_i = 1,\quad \sum_{ij}b_ia_{ij} = \frac{1}{2},\quad\sum_{ijk}b_ia_{ij}a_{jk} = \frac{1}{6},\quad\sum_{ijk}b_ia_{ij}a_{ik} = \frac{1}{3}.
\end{equation*}
These equations may be shown to be  mutually independent \cite{butcherbook}.
\subsection{Composing Runge-Kutta methods. Processing} Butcher \cite{butcheralgebra} developed an algebraic theory of RK methods. If $\psi_h^{(1)}$ and $\psi_h^{(2)}$ represent two RK schemes with $s_1$ and $s_2$ stages respectively, the composition of mappings $\psi_h^{(2)}\circ \psi_h^{(1)}$ corresponds to a single step of size $2h$ of a third RK scheme (with $s_1+s_2$ stages). As distinct from the $f$-dependent mapping $\psi_h^{(2)}\circ \psi_h^{(1)}$,  this new scheme is completely independent of the system (\ref{eq:ode}) under consideration and, accordingly, it is possible to define an operation of composition between the RK schemes themselves. The study of this binary operation  is facilitated by  considering {\em classes} of equivalent RK schemes rather than the schemes themselves, where a class consists of  a scheme and all others that generate the same $\psi_h$ (note that, for instance,  reordering the stages of the method changes the value of the $a_{ij}$ and $b_i$ but does not lead to an essentially different computation, see \cite{butcherbook}). Furthermore, if $y =\psi_h(x)$ is an RK map, then $x$ may be recovered from $y$ by taking a step of size $-h$ from $y$ with  {\em another} RK scheme, the so-called adjoint of the original scheme \cite{control}. This should be compared with the situation for the flow, where $y =\phi_h(x)$ implies $x = \phi_{-h}(x)$.

In order to see that these developments are of  relevance to practical computation, we discuss briefly the idea of processing also introduced  by Butcher \cite{butcher69}.
If $\chi_h$ is  a near-identity mapping in $\R^D$ and $\psi_h$ represents any one-step integrator for (\ref{eq:ode}), the mapping
\begin{equation*}
\widehat{\psi}_h = \chi_h^{-1}\circ \psi_h\circ \chi_h
\end{equation*}
defines a {\em processed} numerical integrator (this is easily interpreted in terms of a change of variables, see \cite{cheap}). For $m\geq 1$,
$$
\widehat{\psi}_h^m = \big(\chi_h^{-1}\circ \psi_h\circ \chi_h\big)^m = \chi_h^{-1}\circ \psi_h^m\circ \chi_h;
$$
therefore to advance $m$ steps with the method $\widehat{\psi}_h$ one preprocesses the initial condition to find $\chi_h(x_0)$, advances $m$ steps with the original method and then postprocess the numerical solution by applying $\chi^{-1}_h$. Postprocessing is only performed when output is desired, not at every time step. If both $\psi_h$ and $\chi_h$ correspond to RK methods, then so does $\widehat{\psi}_h$  as discussed above. The idea of processing is useful in different scenarios. If $\chi_h$ may be chosen in such a way that $\widehat{\psi}_h$ is more accurate  than the original $\psi_h$, one  obtains extra accuracy at the (hopefully small) price of having to carry out the processing;
$\psi_h$ is said to have {\em effective order} $\widehat\nu$ \cite{bss}
when
$\widehat{\psi}_h$ has order $\widehat{\nu}$ larger than the order $\nu$ of $\psi_h$. As a second possibility, the processed integrator may possess, when $\psi_h$ does not, some of the valuable geometric properties presented below.

\subsection{B-series}Butcher series (B-series for short) were introduced in \cite{HW} as a means to systematize the derivation and use of expansions like (\ref{eq:psi}) or (\ref{eq:phi}).
It is convenient to introduce an empty rooted tree $\emptyset$ with  $\sigma(\emptyset) = 1$, $\emptyset ! = 1$, $\eldiff_\emptyset(x) = x$ and,  for each RK method, $c_\emptyset = 1$. If $\trees$ represents the set of all rooted trees (including $\emptyset$), then (\ref{eq:psi}) and (\ref{eq:phi}) become respectively:
\begin{equation*}
\psi_h(x) = \sum_{u\in\trees} h^{|u|}\frac{1}{\sigma(u)}\, c_u\, \eldiff_u(x),\qquad
\phi_h(x) = \sum_{u\in\trees} h^{|u|}\frac{1}{\sigma(u)}\, \frac{1}{u!}\, \eldiff_u(x)
\end{equation*}
($|u|$ is the number of vertices of $u$; $|\emptyset|= 0$). The right-hand sides of these equalities provide examples of {\em B-series}: if $\delta \in \R^\trees$ (i.e.\ $\delta$ is a mapping that associates with each
$u\in\trees$  a real number $\delta_u$), the corresponding B-series is, by definition, the formal series\footnote{Rather than using the normalizing factor $1/\sigma(u)$, the original paper \cite{HW} uses an alternative factor. The present normalization simplifies many formulas. }
\begin{equation}\label{eq:Bseries}
B_\delta(x) = \sum_{u\in\trees} h^{|u|}\frac{1}{\sigma(u)}\, \delta_u\: \eldiff_u(x).
\end{equation}
Note that B-series are relative to the system (\ref{eq:ode}) being studied because the elementary differentials change with $f$.

Obviously the set of all B-series is a vector space. A more important algebraic feature is that, if $\gamma\in\R^\trees$ is a family of coefficients with $\gamma_\emptyset = 1$ and $\delta\in\R^\trees$, then the {\em composition} $B_\delta\big(B_\gamma(x)\big)$ is again a B-series $B_\zeta(x)$; furthermore,
the elements $\zeta_u$ of the family $\zeta$ are functions of the elements of $\delta$ and $\gamma$ and do not vary with $h$ or $f$. For instance for the rooted tree $\ts{\rtg}$,
\begin{equation}\label{eq:decomp}
\zeta_{\ts[0.2]{\rtg}} = \delta_\emptyset\: \gamma_{\ts[0.2]{\rtg}} + \delta_{\ts{\rta}}\: \gamma_{\ts{\rta}}\: \gamma_{\ts{\rtb}}+
 \delta_{\ts{\rtb}}\: \gamma_{\ts{\rtb}}+ \delta_{\ts{\rtb}}\: \gamma_{\ts{\rta}}\: \gamma_{\ts{\rta}}+
 \delta_{\ts{\rtc}}\: \gamma_{\ts{\rta}}+
 \delta_{\ts{\rtd}}\: \gamma_{\ts{\rta}}
 +\delta_{\ts[0.2]{\rtg}}\: \gamma_\emptyset.
\end{equation}
The rooted tree in the left-hand side has been \lq pruned\rq\ in all possible ways (including  complete uprooting $\ts[0]{\rtg}\rightarrow \emptyset$ and no pruning $\ts[0]{\rtg}\rightarrow \ts[0]{\rtg}$); in the right-hand side $\delta$ is evaluated at the rooted tree that remains after pruning and $\gamma$ is evaluated at the pieces that have been removed.

As a first example of the use of B-series, we outline the derivation of the RK expansion (\ref{eq:psi}). We begin by assuming  that the
RK solution $x_{n+1}$ and each internal stage $X_{n,i}$ are expressed as B-series (evaluated at $x_n$) with undetermined coefficients $c\in\R^\trees$, $C^i\in\R^\trees$. Those series are then substituted in (\ref{eq:rk}) and the composition rule we have just described is used to find the B-series for  $f(X_{n,i})$; this procedure yields equations that determine the elementary weights $c_u$. The expansion  (\ref{eq:phi}) may be derived analogously.\medskip

\subsection{The Butcher group}
If $\gamma_\emptyset = 1$ and $B_\delta\big(B_\gamma(x)\big) = B_\zeta(x)$ as above,  we may write $\zeta =\delta\star\gamma$ (recall $\zeta$ does not change with $f$ or $h$). The set of all $\gamma\in\R^\trees$ with $\gamma_\emptyset = 1$ is a group $\G_B$ for the operation $\star$: the so-called {\em Butcher group}. The (Taylor expansion of the) true solution flow $\phi_h$ and (the Taylor expansion of) the RK map $\psi_h$ correspond to the elements of $\G_B$  $\{1/u!\}$ and $\{c_u\}$ respectively; equivalent RK schemes give rise to the same element of $\G_B$. For points in $\G_B$ associated with RK schemes the operation $\star$ obviously corresponds to the composition of RK schemes defined by Butcher.

Assume that the system (\ref{eq:ode}) undergoes an affine change of variables $x = M\overline x+c$ and becomes $(d/dt) \overline x = \overline f(\overline x)$, with $\overline f(\overline x) = M^{-1}f(Mx+c)$. If $\gamma \in\G_B$, $B_\gamma(x) = M \overline B_\gamma(\overline x)+c$, where the
notation $\overline B$ indicates that the elementary differentials in the series are based on $\overline f$; in this way $B_\gamma$ is affine equivariant. In general,\footnote{Of course the B-series for $\phi_h$ is equivariant with respect to all changes of variables.} $B_\gamma$ is not equivariant with respect to more general changes $x = \chi(\overline x)$.

Taylor series methods,  multiderivative RK methods and other one-step integrators $x_{n+1} = \psi_h(x_n)$ not of the form (\ref{eq:rk}) are also represented by elements of $\G_B$. The characterization of the family of integrators that correspond to elements of $\G_B$ has been obtained only recently \cite{charact}; note that those integrators necessarily must be affine equivariant.

 A forest is a formal juxtaposition of nonempty rooted trees, such as \ $\{ \ts[0]{\rta} ,\ts[0]{\rtb}, \ts[0]{\rta}, \ts[0]{\rta}\}$ or $\{ \ts[0]{\rtg}, \ts[0]{\rtc}, \ts[0]{\rtb}\}$; it is understood that here the order is irrelevant (e.g. $\{ \ts[0]{\rta}, ,\ts[0]{\rtb}, \ts[0]{\rta}, \ts[0]{\rta}\}$, and $\{ \ts[0]{\rta}, ,\ts[0]{\rta}, \ts[0]{\rta}, \ts[0]{\rtb}\}$, are the same forest). Also considered is an empty forest that contains no rooted trees. Thus the set of forests $\mathcal F$ is the commutative monoid generated by the set of nonempty rooted trees; by taking all formal linear combinations with real coefficients of forests we obtain a commutative, associative algebra $\mathcal A$. The decomposition of each rooted tree that is necessary to compute $B_\delta(B_\gamma(x))$ (for instance from (\ref{eq:decomp}) \ts[0]\rtg\ is decomposed as
 \begin{equation*}
 \emptyset\otimes \ts[0]{\rtg} +
 \ts[0]{\rta}\otimes \ts[0]{\rta}\ts[0]{\rtb}+
 \ts[0]{\rtb}\otimes \ts[0]{\rtb}+
 \ts[0]{\rtb}\otimes \ts[0]{\rta}\ts{\rta}+
 \ts[0]{\rtc}\otimes \ts[0]{\rta}+
 \ts[0]{\rtd}\otimes \ts[0]{\rta}+
 \ts[0]{\rtg}\otimes \emptyset
 )
 \end{equation*}
 defines a coproduct that turns  $\mathcal A$ into a Hopf algebra: the Connes-Kreimer algebra that appears in renormalization in quantum field theory and elsewhere. Since $\mathcal F$ is a basis for the vector space $\mathcal A$, each linear form on $\mathcal A$ may be identified with a mapping $\in\R^{\mathcal F}$. Each element $\gamma\in\G_B\subset \R^{\trees}$ may be extended to a mapping  $\in\R^{\mathcal F}$ by defining $\gamma$ at a forest as $\prod_u \gamma_u$, where the product is extended to all the trees in the forest. In this way the elements of $\G_B$ may be seen as the characters (multiplicative morphisms) of the Hopf algebra $\mathcal A$. The operation $\star$ defined above is the restriction to
 $\G_B$ of the convolution product in the dual of $\mathcal A$.
 In this way Butcher's work anticipated many results that were to be later used in different subfields of mathematics \cite{brouder}.

 \subsection{Modified equations} Since the properties of maps (discrete dynamical systems) are often more difficult to investigate than those of differential equations (continuous dynamical systems),  it may make sense, given a one-step integrator $\psi_h$, to seek a differential system
$(d/dt) x = \widetilde{f}_h(x)$ whose $h$ flow $\widetilde{\phi}_h$ coincides with the map $\psi_h$. While it is well known  \cite{neish} that it is not possible in general to find such $\widetilde{f}_h$,  for typical integrators and smooth $f$, one may construct a formal series
\begin{equation*}
\widetilde{f}_h(x) = \widetilde{f}^{(0)}(x)+ h \widetilde{f}^{(1)}(x)+h^2\widetilde{f}^{(2)}(x)+\cdots,
\end{equation*}
whose formal $h$-flow exactly matches the Taylor expansion of $\psi_h(x)$. For instance, for Euler's rule with $\psi_h(x) = x+hf(x)$,
$\widetilde{f}_h(x)$ is found to be
\begin{equation}\label{eq:modeuler}
 f(x) -\frac{h}{2} \partial_x f(x) [f(x)] +\frac{h^2}{3} \partial_xf(x) [\partial_x f(x)[f(x)]]
+\frac{h^2}{12} \partial_{xx}f(x)[f(x),f(x)]+\cdots.
\end{equation}
The fact that $\widetilde{f}_h(x)-f(x)=\mathcal{O}(h)$ indicates that the integrator has order 1; for an integrator of order $\nu$,
$\widetilde{f}_h(x)-f(x)=\mathcal{O}(h^\nu)$. When the terms of order $h^\mu$, $\mu = 1, 2,\dots$, and higher are suppressed from (\ref{eq:modeuler}), one obtains a vector field $\widetilde{f}_h^{[\mu]}(x)$ whose $h$-flow differs from $\psi_h$ in $\mathcal{O}(h^{\mu+1})$. Thus Euler's rule applied to (\ref{eq:ode}) is an approximation of order $\mu$ to the {\em modified} differential system $(d/dt)x = \widetilde{f}_h^{[\mu]}(x)$; for $\mu$ large and $|h|$ small the modified system may be expected to provide a very accurate description of the behaviour of the numerical solution.  The idea of using modified systems is very old, see e.g.\ the references in \cite{gss}, and, as we shall see later, has gained prominence with the growing interest in geometric integration.

Note that, for $\widetilde{f}_h(x)$ in (\ref{eq:modeuler}), $h\widetilde{f}_h(x)$ is a B-series. In fact \cite{hairer}, \cite{aust}, each $B_\gamma(x)$ with $\gamma\in\G_B$ coincides with the $h$-flow of a uniquely defined vector field $\widetilde{f}_h(x)$ where $h\widetilde{f}_h(x)$ is a B-series $B_{\beta}(x)$ with $\beta_\emptyset =0$. If we denote by $\g_B$ the set of such $\beta$'s, the mapping $\gamma\mapsto\beta$ is a bijection from $\G_B$ onto $\g_B$. Thus B-series integrators may be handled either by using the coefficients $\{\gamma_u\}$ of the B-series for $\psi_h(x)$ or by the coefficients $\{\beta_u\}$ for the corresponding $h\widetilde{f}_h(x)$. The second option is often advantageous because, while $\g$ is a linear space, $\G_B$ is not; this implies that in many situations properties that are nonlinear when expressed in terms of the $\gamma_u$ become linear for the corresponding $\beta_u$. Formally $\G_B$ is a Lie group and $\g_B$ is its Lie algebra. The Lie bracket in $\g_B$ maybe expressed in terms of the convolution product $\star$ in the dual of $\mathcal A$ (similar developments for the case of words will be presented in the next section).

\subsection{The Hamiltonian case. Geometric integration}In many applications, the system (\ref{eq:ode}) is Hamiltonian, i.e.\ the dimension $D$ is even and the vector field $f(x)$ is of the form $J^{-1}\nabla H(x)$, where $H$ is the (real-valued) Hamiltonian function and $J$ is the matrix
\begin{equation*}
J = \left[ \begin{matrix} O&I\\ -I&O\end{matrix}\right]
\end{equation*}
(the four blocks are of size $D/2\times D/2$). Hamiltonian systems are characterized by the property  that their solution flow $\phi_h$ is, for each $h$, a {\em canonical} or {\em symplectic} transformation, i.e.  (at each $x$)
\begin{equation*}
\big(\partial_x \phi(x)\big)^T J \:\partial_x \phi(x) =J.
\end{equation*}
It is useful, particularly in the context of long-time simulations, to consider one-step integrators $x_{n+1} =\psi_h(x_n)$ that when applied to Hamiltonian system originate a  transformation
$\psi_h$ that is likewise canonical. These integrators are called symplectic \cite{ssc}, \cite{hlw}, \cite{feng} and their importance was first highlighted by  Feng Kang. It was proved independently in \cite{lasagni}, \cite{bit}, \cite{suris} that the RK method (\ref{eq:rk}) is symplectic if
\begin{equation*}
b_ia_{ij}+b_ja_{ji} = b_i b_j,\qquad 1\leq i,j\leq s.
\end{equation*}

In the spirit of the material above, it is sometimes useful to check symplecticness by looking at the B-series of the method rather than by examining the method coefficients. In \cite{canonical} it was proved that, for Hamiltonian systems (\ref{eq:ode}), the  B-series (\ref{eq:Bseries}) with $\delta_\emptyset = 1$ (not necessarily associated with an RK integrator) is symplectic if and only if, for each pair of nonempty rooted trees $u$, $v$,
\begin{equation}\label{eq:sympgroup}
\delta_{u\circ v}+\delta_{v\circ u} = \delta_u\delta_v.
\end{equation}
Here $u\circ v$ denotes the so-called Butcher product of $u$ and $v$, i.e.\ the rooted tree obtained by grafting the root of $v$ into the root of $u$
(e.g.\  $\ts[0]{\rta}\circ \ts[0]{\rtb} = \ts[0]{\rtc}$, $\ts[0]{\rtb}\circ \ts[0]{\rta} = \ts[0]{\rtd}$).  In \cite{hmss} the characterization
(\ref{eq:sympgroup}) was used as a stepping stone to establish the nonexistence of symplectic multiderivative RK schemes.

In terms of  modified vector fields,  symplectic integrators $\psi_h$ may be characterized as methods that when applied to a Hamiltonian problem
(\ref{eq:ode}) result in a  Hamiltonian $\widetilde{f}_h(x)$. Thus symplectic integrators may be alternatively seen as those
integrators that, when applied to  a Hamiltonian system, generate a map $\psi_h$ that formally coincides with the flow of {\em another} Hamiltonian system, hopefully close to the true Hamiltonian. For nonsymplectic schemes $\psi_h$ will coincide with the flow of a vector field perhaps close to $f$ but not within the Hamiltonian class. This interpretation is crucial in symplectic integration \cite{ssc}, \cite{hlw}.

The set of all B-series  $h\widetilde{f}_h(x)$ that are Hamiltonian when $f$ is Hamiltonian provides a Lie subalgebra $\g_0$ of $\g_B$; the associated Lie subgroup of $\G_B$ corresponds to of all symplectic $\psi_h$. An element $\beta\in \g_B$ belongs to $\g_0$ if and only if,
for nonempty $u$ and $v$,
\begin{equation*}
\beta_{u\circ v}+\beta_{v\circ u} = 0.
\end{equation*}
Note that this relation for the Lie algebra is linear, as distinct from the corresponding relation (\ref{eq:sympgroup}) for the group.
For $\beta\in\g_0$ the Hamiltonian function $\widetilde{H}_h(x)$ for the Hamiltonian vector field $\widetilde{f}_h(x)$ is given by a formal series somehow similar to (\ref{eq:Bseries}) but based on so-called (scalar) elementary Hamiltonians \cite{hairer} rather than on (vector-valued) elementary differentials. While there is an elementary differential per rooted tree, there are \lq fewer\rq\ elementary Hamiltonians (one per so-called nonsuperfluous free  tree  \cite{abia}); this explains that for symplectic RK schemes the order conditions corresponding to different rooted trees are not independent.

The study of symplectic integrators for Hamiltonian problems was the first step in what was  termed in \cite{gi} {\em geometric integration}: the integration of differential equations by schemes that preserve important geometric features of the system being integrated \cite{hlw}.
In this way, the literature has envisaged volume preserving integrators to integrate divergence free differential equations,  integrators that preserve relevant invariants, etc. Formal series have been a key element in these studies, see e.g.\
 \cite{cell}, \cite{cfm}, \cite{chartmur}.\medskip

\subsection{Extensions}There are several useful extensions of the Runge-Kutta format (\ref{eq:rk}). In some problems the components of the vector $x$ come to us divided into two or more different groups (for instance in mechanical problems we may have positions and velocities). In those circumstances we may use different coefficients $b_i$, $a_{ij}$ for different groups of components; the result is a {\em partitioned} RK scheme. An important particular case arises when a second-order system $(d^2/dt^2) y = F(y, (d/dt)y)$ is rewritten as a first-order system for $x = (y, (d/dt)y)$; this is the realm of {\em Runge-Kutta-Nystr\"{o}m} methods. In other instances the vector field in (\ref{eq:ode}) may be decomposed as a sum of $N$ parts
\begin{equation*}
f(x) = f_1(x)+\cdots+f_N(x)
\end{equation*}
and we may resort to additive RK methods \cite{aderito} based on evaluations of the individual parts $f_i$ rather than on evaluations of $f$.

The material outlined in previous paragraphs may be adapted to cover all those extensions. A unifying technique has been given in \cite{tesis}. Let us take the additive case as an example.  When Taylor expanding the map $\psi_h$, we find elementary differentials like $\partial_x f_1(x) [f_2(x)]$ or $\partial_{xx} f_2(x) [f_1(x), f_2(x)]$. The structure of these is captured by
{\em coloured} rooted trees, i.e.\ rooted trees where each vertex has been marked  with one of the symbols (colours) $1$, \dots, $N$ (for the elementary differential $\partial_{xx} f_2(x) [f_1(x), f_2(x)]$, the root of \ts[0]{\rtd} is coloured as $2$ and the two terminal vertices as $1$ and $2$). In the definition  (\ref{eq:Bseries}) one has to replace
$\trees$ by the set of all coloured rooted trees; after that, all the developments above are easily adapted to the additive case \cite{aderito}.

In addition to the composition law for B-series that has been the key element above, a {\em substitution law} has been introduced \cite{gilles}, \cite{kur}.

The use of formal series based on rooted trees goes beyond integrators based on the RK idea of repeated evaluations of the vector field. Thus the paper \cite{phil} studies order conditions for splitting and composition methods.

To conclude this section we mention the possibility of using the rooted tree machinery to perform efficiently high-order averaging in periodically or quasiperiodically forced dynamical systems \cite{part1}, \cite{part2}, \cite{orlando}, \cite{part3}. The idea is to expand the  solution as a B-series with oscillatory coeffients and show that those coefficients may be interpolated by nonoscillatory functions of $t$. As a byproduct one may obtain in some circumstances formal conserved quantities.
These are other applications to nonnumerical mathematics of the series methodologies developed to analyze numerical integrators.

\section{Word series}

While the material above dealt with series paratemerized by rooted trees or colour\-ed rooted trees, we now turn our attention to series parameterized by words from an alphabet.

\subsection{Definition}In many situations (see \cite{anderfocm}, \cite{words}) the problem to be integrated is of the form
\begin{equation}\label{eq:odewords}
\frac{d}{dt} x= \sum_{a\in A} \lambda_a(t) f_a(x),\qquad x(0) = x_0,
\end{equation}
where $A$ is a  finite or infinite countable set of indices and, for each $a\in A$, $\lambda_a$ is a scalar-valued function  and $f_a$  a  $D$-vector-valued map. The simplest  example is furnished by the
   autonomous system
\begin{equation}\label{eq:odewordsexample}
\frac{d}{dt} x= f_a(x)+f_b(x),
\end{equation}
where $A =\{a,b\}$ and
$\lambda_a(t) = \lambda_b(t)= 1$. A more complicated nonautonomous example will appear in the next section.

The solution of (\ref{eq:odewords}) has the formal expansion \cite{words}:
\begin{equation}\label{eq:expan}
 x(t) = x_0 + \sum_{n=1}^\infty \sum_{a_1,\dots,a_n\in A}\alpha_{a_1\cdots a_n}(t) f_{a_1\cdots a_n}\!(x_0),
\end{equation}
where  the  mappings
$ f_{a_1\cdots a_n}\!(x)$ and the scalar  functions $\alpha_{a_1\cdots a_n}$ are defined by  the recursions
\begin{equation}\label{eq:wbf}
f_{a_1\cdots a_n}\!(x) = \partial_x f_{a_2\cdots a_n}\!(x)\,f_{a_1}\!(x), \quad n >1,
\end{equation}
and
\begin{equation}\begin{split}
\alpha_{a_1}\!(t) & = \int_0^t \lambda_{a_1}\!(t_1)\,dt_1,\\
\alpha_{a_1\cdots a_n}\!(t) & = \int_0^t \lambda_{a_n}\!(t_n)\,\lambda_{a_1\cdots a_{n-1}}\!(t_n)\,dt_n, \quad n >1.
\end{split}
\label{eq:alfa}
\end{equation}

For the particular case (\ref{eq:odewordsexample}),
 the inner sum in (\ref{eq:expan}) comprises $2^n$ terms and each of them has a coefficient $t^n/n!$: the expansion (\ref{eq:expan}) is  the Taylor series for $x(t)$ as a function of $t$ written in terms of the pieces $f_a$ and $f_b$ rather than in terms of $f$. If the flows $\phi_h^a$ and $\phi_h^b$ of the split systems $(d/dt)x = f_a(x)$ and $(d/dt)x = f_b(x)$ are available analytically (or may be easily approximated numerically), it is often advantageous to consider integrating (\ref{eq:odewordsexample}) by means of {\em splitting} integrators of the form
\begin{equation}\label{eq:splittingalg}
 \psi_h = \phi^b_{d_sh}\circ \phi^a_{c_sh}\circ \cdots \circ\phi^b_{d_1h}\circ \phi^a_{c_1h}
\end{equation}
with $c_i$, $d_i$, $1\leq i\leq s$, constants that specify the method. The Lie-Trotter $\phi^b_h\circ\phi^a_h$ and Strang $\phi^a_{h/2}\circ\phi^b_h\circ\phi^a_{h/2}$ splittings provide the simplest examples. The Taylor expansion of $\psi_h$ involves the individual pieces $f_a$ and $f_b$ and it then makes sense to write the Taylor expansion of the solution flow in the form (\ref{eq:expan}) we have just considered. Of course splitting integrators
that use the flow of the vector fields $f_a$ and $f_b$ are not to be confused with additive RK schemes that just avail themselves of the capability of evaluating the fields $f_a$, $f_b$.

The notation in (\ref{eq:expan}) may be made slightly more compact by considering $A$ as an {\em alphabet}
and the strings $a_1\cdots a_n$ as {\em words}. Then, if $\W_n$ represents the set of all words with
$n$ letters, (\ref{eq:expan}) reads
\[
 x(t) = x_0 + \sum_{n=1}^\infty \sum_{w\in \W_n}\alpha_w\!(t)\, f_w\!(x_0).
\]
If we furthermore introduce the empty word $\emptyset$ and set $\W_0 = \{\emptyset\}$, $\alpha_\emptyset =1$, $f_\emptyset(x) =x$, then
the last expansion becomes
\begin{equation}\label{eq:expanw}
 x(t) = \sum_{n=0}^\infty \sum_{w\in \W_n}\alpha_w\!(t)\, f_w\!(x_0) = \sum_{w\in \W}\alpha_w\!(t)\, f_w\!(x_0),
\end{equation}
where $\W$ represents the set of all words. This suggests   the following definition: If $\delta$ maps $\W$ into $\C$ (i.e.\ $\delta\in\C^\W$),\footnote{While in the case of B-series we only considered real-valued families of coefficients $\delta$, it will be convenient later to consider  word series with complex coefficients.} then its
corresponding {\em word series} is the formal series
\begin{equation}\label{eq:ws}
W_\delta(x) = \sum_{w\in\W} \delta_w f_w\!(x).
\end{equation}
The scalars $\delta_w$ and the functions $f_w$ will be called the {\em coefficients} of the series and the {\em word-basis functions} respectively.
Thus, for each fixed $t$, (\ref{eq:expanw}) is the word series with coefficients $\alpha_w(t)$. As we shall see below, in the particular case (\ref{eq:odewordsexample}), the mapping $\psi_h$ that represents the splitting method (\ref{eq:splittingalg}) corresponds, for each fixed $h$, to a word series.

The definition of word series is clearly patterned after that of B-series. Note however that in (\ref{eq:ws}) the coefficients $\delta_w$ play the role that in
(\ref{eq:Bseries}) is played by $h^{|u|} \delta_u$. As we have pointed out already, the expansion (\ref{eq:expan}) corresponds to a {\em family} of word series: one for each fixed value of $t$. The reason for this small lack of parallelism between the definition of B- and word series is that
while (\ref{eq:psi}) or (\ref{eq:phi}) depend on $h$ through powers $h^j$ which may be made to feature in the definition,
the time dependence of (\ref{eq:expanw}) and related expansions changes with the functions $\lambda_a(t)$.

Each  word-basis function $f_w$, $w\neq \emptyset$, is build up from partial derivatives of the $f_a$, $a\in A$, e.g., if $a,b,c \in A$,
\begin{equation*}
\begin{split}
f_{ba}(x) &= \partial_xf_a(x)\, f_b(x),\\
f_{cba}(x) &= \partial_xf_{ba}(x)f_c(x) = \partial_{xx} f_a(x)[f_b(x),f_c(x)]+\partial_xf_a(x)\,\partial_xf_b(x)\, f_c(x).
\end{split}
\end{equation*}
 Clearly, $\partial_xf_a(x)\, f_b(x)$, $\partial_{xx} f_a(x)[f_b(x),f_c(x)]$, $\partial_xf_a(x)\,\partial_xf_b(x)\, f_c(x)$  are elementary differentials based on rooted trees coloured by the letters of $A$. Thus by expanding each word-basis function in terms of  elementary differentials it would be possible in principle to avoid the introduction of words and  work with coloured rooted trees. However such a move would not be always be advisable
because word series are more {\em compact} than B-series (there are \lq fewer\rq\ word basis functions than elementary differentials) and, additionally, the composition rule for word series is simpler than its counterpart for B-series.

\subsection{The word series group $\G$}

Given $\delta,\delta^\prime\in\C^\W$, we  associate with them its {\em convolution product} $\delta\star\delta^\prime\in\C^\W$ defined by
\begin{equation}\label{eq:convol}
(\delta\star\delta^\prime)_{a_1\cdots a_n} = \delta_\emptyset\delta^\prime_{a_1\cdots a_n}
+ \sum_{j=1}^{n-1} \delta_{a_1\cdots a_j}\delta^\prime_{a_{j+1}\cdots a_n}
+\delta_{a_1\cdots a_n}\delta^\prime_\emptyset
\end{equation}
(here it is understood that $(\delta\star\delta^\prime)_\emptyset = \delta_\emptyset\delta^\prime_\emptyset$). The convolution product is not commutative, but it is associative and has a unit (the element $\uno\in \C^\W$ with $\uno_\emptyset = 1$ and $\uno_w = 0$ for $w\neq \emptyset$).

If $w\in\W_m$ and $w^\prime\in \W_n$ are words, $m,n\geq 1$, its {\em shuffle product} $w\sh w^\prime$ \cite{reu} is the formal sum of the $(m+n)!/(m!n!)$
words with $m+n$ letters that may be obtained by interleaving the letters of $w$ and $w^\prime$ while preserving the order in which the letters appear in each word.
In addition $\emptyset \sh w = w  \sh \emptyset= w$ for each $w\in\W$. The operation $\sh$ is commutative and associative and has the word $\emptyset$ as a unit. We  denote by $\G$ the set of those $\gamma \in \C^\W$  that satisfy the so-called {\em shuffle relations:} $\gamma_\emptyset = 1$ and, for each $w,w^\prime\in \W$,
\begin{equation*}
\gamma_w\gamma_{w^\prime} = \sum_{j=1}^N \gamma_{w_j}\qquad \mbox{\rm if}\qquad w\sh w^\prime = \sum_{j=1}^N w_j.
\end{equation*}
The set $\G$ with the operation $\star$ is a (non-commutative) formal Lie group, which plays here the role played by $\G_B$ in the preceding section. For each fixed $t$, the family of coefficients defined by (\ref{eq:alfa}) and $\alpha_\emptyset(t) =1$ is an element of the group $\G$, \cite[Corollary 3.5]{reu}. As we shall see below, some numerical integrators for (\ref{eq:odewords}) are also associated with families of elements of $\G$ parameterized by the stepsize.

For $\gamma \in \G$, the word series $W_\gamma(x)$ is equivariant with respect to {\em arbitrary } (not necessarily affine) changes of variables $x=\chi(\overline x)$ \cite[Proposition 3.1]{part2}. In fact, if word series are rewritten as B-series (with colours from $A$), then the family of word series $W_\gamma(x)$ with $\gamma\in \G$ exactly corresponds to the family of B-series that are equivariant with respect to arbitrary changes of variables. By implication only integrators that are equivariant with respect to arbitrary changes of variables are candidates to possess a word-series expansion with coefficients in $\G$. Thus $\G$ may be regarded as a (small) subgroup of $\G_B$.

In analogy with the composition of B-series, for $\gamma\in\G$, $W_\gamma(x)$ may be substituted in an arbitrary word series $W_\delta(x)$, $\delta\in\C^\W$, to get a new word series; more precisely
\begin{equation}\label{eq:act}
W_\delta\big (W_{\gamma}(x)\big) = W_{\gamma\star \delta}(x),
\end{equation}
i.e.\ the coefficients of the word series resulting from the substitution are given by the convolution product $\gamma\star\delta$. It follows immediately that, in the particular case  (\ref{eq:odewordsexample}), the map $\psi_h$ defined in (\ref{eq:splittingalg}) corresponds, for each fixed $h$, to a word series with coefficients in $\G$.

The Lie algebra $\g$ of the group $\G$ consists of the elements $\beta\in\C^\W$ such that $\beta_\emptyset = 0$ and, for each pair of nonempty words $w,w^\prime$,
\[
\sum_{j=1}^N \beta_{w_j} = 0\qquad \mbox{\rm if}\qquad w\sh w^\prime = \sum_{j=1}^N w_j.
\]
The bracket operation in $\g$ is
\begin{equation*}
[\beta,\beta^\prime] = \beta\star\beta^\prime-\beta^\prime\star\beta
\end{equation*}
 and corresponds to the Jacobi bracket (commutator) of the associated word series, i.e.\ for $\beta,\beta^\prime\in \g$:
\[
\big(\partial_x W_{\beta^\prime}(x)\big) W_\beta(x) - \big(\partial_x W_\beta(x)\big) W_{\beta^\prime}(x)= W_{[\beta,\beta^\prime]}(x).
\]
Since for $\beta\in\g$, the word series $W_\beta(x)$ belongs to the Lie algebra (for the Jacobi bracket) generated by the  $f_a$, the Dynkin-Specht-Wever formula \cite{jac} may be used to rewrite the word series in terms of iterated commutators of these mappings:
\begin{equation*}
W_\beta(x)
 = \sum_{n=1}^\infty \frac{1}{n} \sum_{a_1,\dots,a_n\in A}
 \beta_{a_1\cdots a_n} [[\cdots[[f_{a_1},f_{a_2}],f_{a_3}]\cdots],f_{a_n}](x).
\end{equation*}
(For $n=1$ the terms in the inner sum are of the form $\beta_{a_1} f_{a_1}(x)$.)

In the particular case where each $f_a$ is a Hamiltonian vector field, the iterated commutators are well known to  be Hamiltonian and therefore so is
$W_\beta(x)$. Furthermore the Hamiltonian function
of the vector field $W_\beta(x)$  is
$$
\Ham_\beta(x) = \sum_{w\in\W,\, w\neq\emptyset} \beta_w H_w(x),
$$
where, for each nonempty word $w=a_1\cdots a_n$,
\begin{equation}\label{eq:wordham}
H_w(x) = \frac{1}{n}\{\{\cdots\{\{H_{a_1},H_{a_2}\},H_{a_3}\}\cdots\},H_{a_n}\}(x).
\end{equation}
Here $\{\cdot,\cdot\}$ is the Poisson bracket  defined by $\{A,B\}(x) = \nabla A(x)^T J^{-1} B(x)$.

We conclude this subsection with an interpretation of the material above in terms of Hopf algebras. The product $\sh$ may be extended in a bilinear way from words to linear combinations of words. After such an extension, the vector space  $\C\langle A\rangle $ of all such  linear combinations   is  a unital, commutative, associative algebra, {\em the shuffle algebra}, denoted by sh$(A)$ (see \cite{reu},  \cite{anderfocm}). Deconcatenation defines a coproduct and turns sh$(A)$ into a (commutative, connected, graded) {\em Hopf algebra} \cite{brouder}. The sets $\G$ and $\g$  are then respectively the  group of characters and the Lie algebra of infinitesimal characters of the Hopf algebra sh$(A)$.

\section{Extended word series}

Extended word series, introduced in \cite{words}, are a generalization of word series  to cope with perturbed integrable problems and its discretizations.

\subsection{Definition} We now consider systems of the form
\begin{equation}\label{eq:odefriday}
\frac{d}{dt} \left[ \begin{matrix}y\\ \theta\end{matrix}\right]
= \left[ \begin{matrix}0\\ \omega\end{matrix}\right]
+f(y,\theta),
\end{equation}
where $y\in\R^{D-d}$, $0<d\leq D$, $\omega\in\R^d$ is a vector of frequencies $\omega_j>0$, $j = 1,\dots, d$, and $\theta$ comprises $d$ angles, so that $f(y,\theta)$ is $2\pi$-periodic in each component of $\theta$ with Fourier expansion
$$
f(y,\theta) = \sum_{\bk \in\Z^d} \exp(i \bk\cdot \theta)\: \hat f_\bk(y)
$$
($\hat f_\bk(y)$ and $\hat f_{-\bk}(y)$ are mutually conjugate, so as to have a real problem). Systems of this form appear in many applications, perhaps after a change of variables. For instance, any system  $(d/dt)z = Mz+F(z)$, where
$M$ is a skew-symmetric constant matrix, may be brought to the format (\ref{eq:odefriday}); other examples are discussed in \cite{words}. When $f\equiv 0$ the system is integrable (the angles rotate with uniform angular velocity and $y$ remains constant) and accordingly we refer to problems of the form (\ref{eq:odefriday}) as perturbed integrable problems and to $f$ as the perturbation (some readers may prefer to substitute $\epsilon f$ for $f$).

After introducing the  functions
\begin{equation}\label{eq:fbk}
f_\bk(y,\theta) = \exp(i\bk\cdot \theta)\: \hat f_\bk(y),\qquad y\in\R^{D-d},\: \theta\in\R^d,
\end{equation}
we have
\begin{equation}\label{eq:ode2}
\frac{d}{dt} \left[ \begin{matrix}y\\ \theta\end{matrix}\right]
= \left[ \begin{matrix}0\\ \omega\end{matrix}\right]
+f(y,\theta)
= \left[ \begin{matrix}0\\ \omega\end{matrix}\right]
+\sum_{\bk \in\Z^d} f_\bk(y,\theta).
\end{equation}

To find the solution with initial conditions
\begin{equation}\label{eq:ic2}
y(0) = y_0,\qquad \theta(0) = \theta_0,
\end{equation}
we perform the time-dependent change of variables $\theta =\eta+t\omega $ to get
\begin{equation*}
\frac{d}{dt} \left[ \begin{matrix}y\\ \eta\end{matrix}\right]
= \sum_{\bk \in\Z^d} \exp(i \bk\cdot \omega t)\:f_\bk(y,\eta),
\end{equation*}
a particular instance of (\ref{eq:odewords}).
 The formula (\ref{eq:expan}) yields
\begin{equation*}
\left[ \begin{matrix}y(t)\\ \eta(t)\end{matrix}\right]
=
\left[ \begin{matrix}y(0)\\ \eta(0)\end{matrix}\right]
+\sum_{n=1}^\infty \sum_{\bk_1,\dots,\bk_n}\alpha_{\bk_1\cdots \bk_n}(t)\, f_{\bk_1\cdots \bk_n}(y(0),\eta(0)),
\end{equation*}
where the coefficients $\alpha$ are still given by (\ref{eq:alfa}) (with $\lambda _\bk(t) = \exp(i \bk\cdot \omega t))$ and the word basis functions are defined by (\ref{eq:fbk}) and (\ref{eq:wbf}) (the Jacobian in (\ref{eq:wbf}) is of course taken with respect to the $D$-dimensional variable $(y,\theta)$). We conclude that, in the original variables, the solution flow of (\ref{eq:ode2}), has the formal expansion
\begin{equation}\label{eq:solucion}
\phi_t(y_0,\theta_0) =
\left[ \begin{matrix}y(t)\\ \theta(t)\end{matrix}\right] =
\left[ \begin{matrix}y_0\\ \theta_0\end{matrix}\right]
+\left[ \begin{matrix}0\\t \omega \end{matrix}\right]
+\sum_{n=1}^\infty \:\sum_{\bk_1,\dots,\bk_n}\alpha_{\bk_1\cdots \bk_n}(t)\, f_{\bk_1\cdots \bk_n}(y_0,\theta_0).
\end{equation}
Note that the word basis functions are {\em independent of the frequencies} $\omega$ and the coefficients $\alpha$ are {\em independent of $f$}. 

With the notation of the preceding section, we write (\ref{eq:solucion}) in the following form (here and later $x = (y,\theta)$):
$$
x(t) = \left[ \begin{matrix}0\\t \omega \end{matrix}\right] +W_{\alpha(t)}(x_0).
$$
 We introduce the vector space  $\calC = \C^d\times \C^\W$ and,
 if $(v,\delta)\in\calC$, define its
corresponding {\em extended word series} to be the formal series

\[
\overline{W}_{(v,\delta)}(x) =  \left[ \begin{matrix}0\\ v \end{matrix}\right] +\sum_{w\in\W} \delta_w f_w\!(x).
\]

Then the solution (\ref{eq:solucion}) of (\ref{eq:ode2})--(\ref{eq:ic2}) has the expansion
$$x(t) = \overline{W}_{(t\omega,\alpha(t))}(x_0),\qquad (t\omega,\alpha(t)) \in \calC,$$
with  $\alpha(t)\in \G \subset \C^\W$  as defined before.

\subsection{The extended word series group $\overline\G$}The symbol $\overline{\G}$ denotes the subset of $\calC$ comprising the elements $(u,\gamma)$ with $u \in \C^d$ and $\gamma\in\G$. For each $t$, the
solution coefficients $(t\omega,\alpha(t))\in\calC$ found above provide an example of element of $\overline{\G}$. Some numerical integrators for
(\ref{eq:odefriday}) will be shown below to have, for each value of the stepsize $h$, an expansion in extended word series with coefficients in $\overline{\G}$.

We introduce linear opertors $\Xi_{v}$, $\xi_{v}$ as follows. If $v$ is a $d$-vector,
 $\Xi_{v}$ is the linear operator in $\C^\W$ that maps each $\delta \in \C^\W$ into the element of $\C^\W$ defined by
 $(\Xi_{v} \delta)_\emptyset = \delta_\emptyset$ and
\begin{equation*}
(\Xi_{v}\delta)_w =\exp( i (\bk_1+\cdots+\bk_n)\cdot v)\: \delta_w,
\end{equation*}
 for $w = \bk_1\dots\bk_n$. For the linear operator $\xi_{v}$ on $\C^{\W}$, $(\xi_{v} \delta)_{\emptyset}=0$, and for each word $w=\bk_1 \cdots \bk_n$,
\begin{equation*}
(\xi_{v} \delta)_{w} = i  (\bk_1 +\cdots + \bk_n)\cdot v\: \delta_{w}.
\end{equation*}

With the help of $\Xi_u$ we define a binary operation $\bigstar$:
 if $(u,\gamma)\in\overline{\G}$ and $(v,\delta) \in\calC$, then
 $$
(u,\gamma) \bigstar(v,\delta) = (\gamma_\emptyset v + \delta_{\emptyset} u,
\gamma \star (\Xi_{u} \delta))\in\calC.
$$
By using (\ref{eq:act}), it is a simple exercise to check that  $\overline{\G}$ acts by {\em composition} on extended word series as follows:
\begin{equation*}
\overline{W}_{(v,\delta)}\big(\overline{W}_{(u,\gamma)}(x)\big) =
\overline{W}_{(u,\gamma)\bigstar (v,\delta)}(x), \qquad \gamma\in{\G}.
\end{equation*}
In fact we have defined the operation $\bigstar$  so as to ensure this property.
The set $\overline{\G}$ is a group for the product $\bigstar$ and $\C^d$ and $\G$ are subgroups of $\overline{\G}$.\footnote{
Consider the group homomorphism  from the additive group $\C^d$ to the group of automorphisms of $\G$ that maps each $\mu \in \C^d$ into $\Xi_{\mu}$. Then $\overline{\G}$ is the (outer) semidirect product of $\G$ and the additive group $\C^d$ with respect to this homomorphism.}   The unit of $\overline{\G}$ is the element $\overline\uno = (0,\uno)$.

As a set, the Lie algebra $\overline{\g}$ of the group $\overline{\G}$ consists of the elements
$(v,\delta)\in\calC$ with $\delta\in\g$.
For $(v,\delta),(u,\eta) \in \C^d\times\g$, the Jacobi bracket of the vector fields $\overline{W}_{(v,\delta)}$, $\overline{W}_{(u,\eta)}$ may be shown to be given by \cite{words}
\begin{equation*}
[\overline{W}_{(v,\delta)},\overline{W}_{(u,\eta)}] = \overline{W}_{(0,\xi_{v}\eta-\xi_{u} \delta + \delta*\eta -\eta*\delta )};
\end{equation*}
accordingly the bracket of $\overline{\g}$ has the expression
\begin{equation*}
[(v,\delta),(u,\eta))] =
(0,\xi_{v}\eta-\xi_{u} \delta + \delta*\eta -\eta*\delta).
\end{equation*}
The $0$ in the right-hand side reflects  the fact that $\C^d$ is an Abelian subgroup of $\overline{\G}$.

Assume that  in (\ref{eq:ode2}) the dimension $D$ is even with  $D/2- d=m\geq0$ and that the vector of variables takes the form
$$
x = (y,\theta) =(p^1,\dots,p^m; q^1,\dots,q^m;a^1,\dots,a^d;\theta^1,\dots,\theta^d),
$$
where $p^j$ is the momentum canonically conjugate to the co-ordinate $q^j$ and $a^j$ is the momentum (action) canonically conjugate to the co-ordinate (angle) $\theta^j$. If each $f_\bk(x)$ in (\ref{eq:fbk}) is a Hamiltonian vector field with Hamiltonian function $H_\bk(x)$, then the system (\ref{eq:ode2}) is itself Hamiltonian for the Hamiltonian function
$$
\sum_{j=1}^d \omega_j a^j +\sum_{\bk\in\Z^d} H_\bk(x).
$$

For each $(\omega,\beta)\in\overline{\g}$, the extended word series $\overline{W}_{(\omega,\beta)}(x)$ is a  Hamiltonian formal vector field, with Hamiltonian function
\begin{equation*}
\sum_{j=1}^d \omega_j a^j + \sum_{w\in\W,\,w\neq\emptyset}\beta_w H_w,
\end{equation*}
with $H_w(x)$ as in (\ref{eq:wordham}).

The paper \cite{words} shows how to use the algebraic rules we have just described to reduce (\ref{eq:odefriday}) to {\em normal} form where all oscillatory components are removed from the solution by means of a suitable change of variables.

 \section{Analyzing splitting methods for perturbed integrable problems}

Splitting algorithms  are natural candidates to integrate (\ref{eq:ode2}).
 Given real coefficients, $a_j$ and $b_j$, $j=1,\dots,r$, we study the splitting integrator
\begin{equation}\label{eq:integrat}
{\psi}_h = \phi^{(P)}_{b_rh}\circ \phi^{(U)}_{a_rh}\circ\cdots\circ \phi^{(P)}_{b_1h}\circ \phi^{(U)}_{a_1h},
\end{equation}
where $\phi_t^{(U)}$ and
$\phi_t^{(P)}$ denote respectively the  $t$-flows of the split systems corresponding to the unperturbed  dynamics
\begin{equation*}
\frac{d}{dt} \left[ \begin{matrix}y\\ \theta\end{matrix}\right]
= \left[ \begin{matrix}0\\ \omega\end{matrix}\right],
\end{equation*}
and the perturbation
\begin{equation*}
\frac{d}{dt} \left[ \begin{matrix}y\\ \theta\end{matrix}\right]
=
f(y,\theta).
\end{equation*}

Since the unperturbed dynamics with frequencies $\omega_j$ is reproduced exactly by (\ref{eq:integrat}), one would naively hope that the accuracy of the integrator would be dictated for the size of $f$ uniformly in $\omega$. It is well known that such an expectation is unjustified, see e.g. \cite{molly1}, \cite{molly2}, since the oscillatory character of the solution leads to a very complex behaviour of the numerical solution. The algebraic machinery of extended word series has been used in \cite{words} to provide a very detailed description of the performance of the integrators; we shall quote below some of the results in  that paper.

\subsection{The truncation error} Clearly, the mapping $\phi_t^{(U)}$ has an expansion in extended word series
$$
\phi_t^{(U)}(x) = \overline{W}_{(t\omega,\uno)}(x),\qquad (t\omega,\uno)\in\overline{\G};
$$
furthermore,
$$
\phi_t^{(P)}(x) = \overline{W}_{(0,\tau(t))}(x),\qquad (0,\tau(t))\in\overline{\G},
$$
where $\tau(t)\in\G$  comprises the  Taylor coefficients, i.e.\ $\tau_w(t)= t^n/n!$ if $w\in\W_n$. It follows that ${\psi}_h$ also possesses an expansion as an extended word series with coefficients in $\overline \G$. A simple computation using the definition of $\bigstar$ shows that:
$${\psi}_h(x) = \overline{W}_{(ha\omega,\widetilde{\alpha}(h))}(x),
$$
where $a = \sum_i a_i$, and
$\widetilde{\alpha}(h) \in\G$ is specified by
$\widetilde{\alpha}_\emptyset(h)=1$ and, for $n= 1,2,\dots$,
\begin{equation*}
\widetilde{\alpha}_{\bk_1 \cdots \bk_n}(h) = h^n\sum_{1\leq j_1 \leq \cdots \leq j_n \leq r} \frac{b_{j_1}\cdots b_{j_n}}{\sigma_{j_1 \cdots j_n}}\,
\exp(i \, (c_{j_1}\bk_1 +\cdots +c_{j_n}\bk_n  ) \cdot \omega h).
\end{equation*}
Here,
 $$c_j = a_1 + \cdots + a_{j}, \qquad 1\leq j\leq r,$$
and,
\begin{eqnarray*}
  \sigma_{j_1\cdots j_n}=\frac{1}{n!} &\mbox{\rm if }& j_1 = \cdots = j_n,\\
  \sigma_{j_1\cdots j_n}=\frac{1}{\ell!}\, \sigma_{j_{\ell+1}\cdots j_n} &\mbox{\rm if }& \ell<n,\quad j_1=\cdots = j_{\ell} < j_{\ell+1} \leq \cdots \leq j_n.
\end{eqnarray*}
By subtracting the extended word series of the true solution and of the integrator, we obtain an expansion of the local error that may be used to investigate the order of consistency \cite{words}. However such an investigation  throws  light on the behaviour of the numerical method only when $|h|$ is small with respect to the periods of the oscillations in the problem. A much better description of the numerical solution may be obtained by means of the modified systems that we describe next.

\subsection{Modified systems}
For $n = 1,2, \dots$, we look for a system
\begin{equation*}
\frac{d}{dt}  x = \overline{W}_{(\omega,\widetilde\beta(h))}( x),\qquad \widetilde\beta(h)\in \g,
\end{equation*}
where $\beta_w(h)= 0$ for words with more than $n$ letters, such that,  for words with $\leq n$ letters, the extended word series expansion of its flow matches that of the integrator.
Such a system may be constructed \cite{words} provided that there is no numerical resonance of order $\leq n$, i.e.\ there is no set $\bk_1$, \dots, $\bk_r$ with $(\bk_1+\cdots +\bk_r)\cdot \omega h=2\pi j$, $r\leq n$,  $j\neq 0$. Furthermore, in the Hamiltonian case, the modified systems will also be Hamiltonian.

In the limit where $n$ increases indefinitely one obtains, if there is no numerical resonance of any order, a modified system whose formal $h$-flow exactly reproduces the extended word series expansion of $\widetilde \phi_h$. As detailed in \cite{words} such modified systems provide a very accurate description of the behaviour of the computed solutions. Among other things, it makes it possible to prove $\mathcal{O}(h)$ error bounds for values of $h$ away from first order numerical resonances but not necessarily small with respect to the periods of the fastest oscillations in the problem.

\section*{Acknowledgements}

A. Murua and J.M. Sanz-Serna have been supported by proj\-ects MTM2013-46553-C3-2-P and MTM2013-46553-C3-1-P from Ministerio de Eco\-nom\'{\i}a y Comercio, Spain. Additionally A. Murua has been partially supported by the Basque Government  (Consolidated Research Group IT649-13).
\section*{References}

\medskip

\renewcommand{\refname}{}    
\vspace*{-36pt}  


\frenchspacing
\begin{thebibliography}{7}

\bibitem{aderito} Araujo, A. L., Murua, A.,   Sanz-Serna, J. M., Symplectic methods based on decompositions, \emph{SIAM J. Numer. Anal.} \textbf{34} (1997), 1926-1947.

\bibitem{brouder}Brouder, Ch., Trees, renormalization and differential equations,\emph{ BIT Numerical Mathematics} \textbf{44} (2004), 425--438.

\bibitem{butcher63}Butcher, J. C.,  Coefficients for the study of Runge-Kutta integration processes, \emph{J. Austral. Math. Soc.} \textbf{3}
            (1963), 185--201.

\bibitem{butcher69}   Butcher, J. C.,  The effective order of Runge-Kutta methods. In
    \emph{Conference on the numerical solution of differential equations} (ed. by J. Ll. Morris). Lecture Notes in Math. Vol. 109, Springer, Berlin, 1969, 133--139.

\bibitem{butcheralgebra} Butcher, J. C., An algebraic theory of integration methods, \emph{Math. Comp.} \textbf{19} (1972), 79--106.

\bibitem{buthis}Butcher, J. C., A history of Runge-Kutta methods, \emph{Appl. Numer. Math.} \textbf{20} (1996), 247--260.

\bibitem{butcherbook} Butcher, J. C., \emph{Numerical Methods for Ordinary Differential Equations, 2nd ed.} Wiley, Chichester, 2008.

\bibitem{bss} Butcher, J. C.,  Sanz-Serna, J. M.,  The number of conditions for a Runge-Kutta method to have effective order $p$, \emph{Appl. Numer. Math.} \textbf{22} (1996), 103--111.

\bibitem{kur}Calaque, D., Ebrahimi-Fard, K., Manchon, D., Two interacting Hopf algebras of trees: A Hopf-algebraic approach to composition and substitution of B-series, \emph{Adv.  Appl. Math.} \textbf{47} (2011), 282--308.

\bibitem{aust}Calvo, M. P., Murua, A., and Sanz-Serna, J. M., Modified equations for ODEs. In Chaotic Numerics (ed. by P. E. Kloeden and K. J. Palmer). Contemporary Mathematics, Vol. 172, American Mathematical Society, Providence, 1944,  63--74.

\bibitem{canonical}  Calvo, M. P., Sanz-Serna, J. M., Canonical B-series, \emph{Numer. Math.} \textbf{67} (1994), 161--175.

\bibitem{cell} Celledoni, E., McLachlan, R. I., Owren, B., and Quispel, G. R. W., Energy-preserving integrators and the structure of B-series, \emph{Found. Comput. Math.} \textbf{10} (2010), 673--693.

\bibitem{cfm} Chartier, P.,  Faou, E.,  Murua, A., An algebraic approach to invariant preserving integrators: The case of quadratic and Hamiltonian invariants, \emph{Numer. Math.}, \textbf{103} (2006), 575--590.

\bibitem{gilles}Chartier, P., Hairer, E., Vilmart, G., Algebraic structures of B-series, \emph{Found. Comput. Math.}\textbf{10}  (2010) 407-–427.

\bibitem{chartmur} Chartier, P., Murua, A., Preserving first integrals and volume forms of additively split systems, \emph{IMA J. Numer. Anal.}
    \textbf{27} (2007), 381--405.

\bibitem{part1}Chartier, P.,  Murua, A., Sanz-Serna, J. M., Higher-Order averaging, formal series and numerical integration I:
     B-series, \emph{Found. Comput. Math.} \textbf{10}  (2010), 695--727.

\bibitem{part2}Chartier, P.,  Murua, A., Sanz-Serna, J. M.,  Higher-Order averaging,
        formal series and numerical integration II: the quasi-periodic
        case,   \emph{Found. Comput. Math.}  \textbf{12}   (2012), 471--508.

\bibitem{orlando}Chartier, P.,  Murua, A., Sanz-Serna, J. M., A formal
    series approach to averaging: exponentially small error estimates,
     \emph{DCDS A}  \textbf{32} (2012), 3009--3027.

\bibitem{part3}Chartier, P.,  Murua, A., Sanz-Serna, J. M.,  Higher-Order averaging,
        formal series and numerical integration III, \emph{Found. Comput. Math.}  (2013) DOI 10.1007/s10208-013-9175-7.

\bibitem{k} Ebrahimi-Fard, K., Lundervold, A.,  Malham, S. J. A., Munte-Kaas, H., Wiese, A.,  Algebraic structure of stochastic expansions and efficient simulation, \emph{Proc. R. Soc. A} \textbf{468} (2012), 2361--2382.

\bibitem{fm} Fauvet, F., and F. Menous, M., Ecalle's arborification-coarborification transforms and Connes-Kreimer Hopf algebra, arXiv; 1212.4740v2.

\bibitem{feng} Feng, H., Qin, M., \emph{Symplectic Geometric Algorithms for Hamiltonian Systems.} Springer, Berlin, 2010.

\bibitem{molly1} Garc\'{\i}a-Archilla, B., Sanz-Serna, J. M.,  Skeel, R. D., Long-time-step methods for oscillatory differential equations, \emph{SIAM J. Sci. Comput.} \textbf{20} (1998), pp.~930--963.

\bibitem{gss}  Griffiths, D. F.,  Sanz-Serna, J. M.,  On the scope of the method of modified equations, \emph{SIAM J. Sci. Statist. Comput.} \textbf{7} (1986), 994--1008.

\bibitem{hairer}Hairer, E., Backward error analysis of numerical integrators and symplectic methods, \emph{Annals Numer. Math.} \textbf{1} (1994), 107--132.

\bibitem{hlw}Hairer, E., Lubich, Ch., Wanner, G., \emph{ Geometric Numerical Integration, 2nd
        ed.} Springer, Berlin, 2006.

\bibitem{hmss} Hairer E., Murua, A., Sanz-Serna, J. M., The nonexistence of symplectic multiderivative Runge-Kutta methods, \emph{BIT}
\textbf{34} (1994), 80--87.

\bibitem{HW} Hairer E.,  Wanner, G.,  On the Butcher group and general multi-value methods, \emph{Computing} \textbf{13} (1974), 1--15.

\bibitem{jac}Jacobson, N., \emph{Lie Algebras.} Dover, New York, 1979.

\bibitem{nuevo} Kawski, M., Sussmann, H. J., Nonommutative power series and formal Lie algebraic techniques in nonlinear control theory. In \emph{ in Operators, Systems, and Linear Algebra}, (ed. by U. Helmke, D. Pratzel-Wolters, E. Zerz). Teubner, Stuttgart, 1997, 111--118.

\bibitem{lasagni}Lasagni, F. M., Canonical Runge-Kutta methods, \emph{ZAMP} \textbf{39} (1988), 952-953.

\bibitem{cheap} Lopez-Marcos, M. A., Skeel, R. D., Sanz-Serna, J. M., Cheap enhancement of symplectic integrators. In \emph{Numerical Analysis 1995} (ed. by D. F. Griffiths and G. A. Watson). Pitman Research Notes in Mathematics 344, Longman Scientific and Technical, London, 1996, 107--122.

\bibitem{charact} McLachlan, R., Modin, K.,
Munthe-Kaas, H., Verdier, O., B-series methods are exactly the local, affine equivariant methods, arXiv.1409v1.

\bibitem{tesis} Murua, A., Formal Series and Numerical integrators. Part I: Systems of ODEs and symplectic Integrators, \emph{Appl. Numer. Math.} \textbf{29} (1999), 221--251.

\bibitem{anderfocm}Murua,A., The Hopf algebra of rooted trees, free Lie algebras and Lie series,
     \emph{Found. Comput. Math.} \textbf{6}  (2006), pp.~387--426.

\bibitem{phil} Murua, A., Sanz-Serna, J. M., Order conditions for numerical integrators obtained by composing simpler integrators, \emph{Phil. Trans. R. Soc. Lond. A} \textbf{357} (1999), 1079--1100.

\bibitem{words} Murua, A., Sanz-Serna, J. M., Word series for dynamical systems and their numerical integrators, arXiv:1502.05528.

\bibitem{neish} Neishtadt, A. I., The separation of motions in systems with rapidly rotating phase, \emph{J. Appl. Math. Mech.} \textbf{48} (1984), 133--139.

\bibitem{reu} Reutenauer, C., \emph{Free Lie Algebras.} Clarendon Press, Oxford, 1993.

\bibitem{bit} Sanz-Serna, J. M., Runge-Kutta schemes for Hamiltonian systems, \emph{BIT} \textbf{28} (1988), 877--883.

\bibitem{gi}  Sanz-Serna, J. M.,  Geometric integration.  In  \emph{The State of the Art in Numerical Analysis} (ed. by I. S. Duff and G. A. Watson). Clarendon Press, Oxford, 1997, 121--143.

\bibitem{molly2} Sanz-Serna, J. M, Mollified impulse methods for highly oscillatory differential equations, \emph{SIAM J. Numer. Anal.}, \textbf{46} (2008), 1040-1059.

\bibitem{control}Sanz-Serna, J. M., Symplectic Runge-Kutta schemes for adjoint equations, automatic differentiation, optimal control and more, arXiv: 1503.04021.

\bibitem{abia}Sanz-Serna, J. M.,  Abia, L., Order conditions for canonical Runge-Kutta schemes, \emph{SIAM J. Numer. Anal.} \textbf{28} (1991), 1081-1096.

\bibitem{ssc}Sanz-Serna, J. M.,   Calvo, M. P., \emph{Numerical Hamiltonian
        problems.} Chapman and Hall, London, 1994.

\bibitem{suris}Suris, Y. B., Preservation of symplectic structure in the numerical solution of Hamiltonian systems. In \emph{Numerical Solution of Differential Equations} (ed. by S. S. Filippov). Akad. Nauk. SSSR, Inst. Prikl. Mat., Moscow, 1988, 138--144 (in Russian).


\end{thebibliography}
\end{document}